\documentclass[11pt]{article}
\usepackage{dsfont, amssymb,amsmath,amscd,latexsym, amsthm, amsxtra,amsfonts}
\usepackage[all]{xy}
\newtheorem{theorem}{Theorem}[section]

\newtheorem{remark}[theorem]{Remark}

\newtheorem{Def}{Definition}[section]
\newtheorem{Them}{Theorem}[section]
\newtheorem{Prop}{Proposition}[section]
\newtheorem{Remark}{Remark}[section]
\newtheorem{Lemma}{Lemma}[section]

\begin{document}
\makeatletter
\def\@setauthors{%
\begingroup
\def\thanks{\protect\thanks@warning}%
\trivlist \centering\footnotesize \@topsep30\p@\relax
\advance\@topsep by -\baselineskip
\item\relax
\author@andify\authors
\def\\{\protect\linebreak}%
{\authors}%
\ifx\@empty\contribs \else ,\penalty-3 \space \@setcontribs
\@closetoccontribs \fi
\endtrivlist
\endgroup } \makeatother
 \baselineskip 16pt
\title{{\Large {\bf Multidimensional SDE with anticipating initial process
and reflection }}\thanks{ This work is
 supported by NSFC and  SRF for ROCS, SEM}}
\author{ Zongxia \ Liang \thanks{ email: zliang@math.tsinghua.edu.cn}
\\ \small Department of Mathematical Sciences,
\  Tsinghua University,\\
\small  Beijing 100084, People's Republic of China }
\date{}
 \maketitle
 \begin{abstract} \noindent In
this paper, the strong solutions $ (X, L)$ of multidimensional
stochastic differential equations with reflecting boundary  and
possible anticipating initial random variables is established. The
key is to obtain some substitution formula for Stratonovich
integrals via a uniform convergence of the corresponding Riemann
sums and to prove continuity of functionals of $ (X, L)$. \\[10pt]
 \end{abstract}
 \setcounter{equation}{0}
\section{ \bf{\small  Introduction and main results}}
Let $\mathcal{O}$ be a  smooth  bounded  open set in $ \Re^d $. $
{\bf n}(x) $ denotes the cone of unit outward normal vectors to $
\partial \mathcal{O} $ at $x$, that is,
\begin{eqnarray}
&(i) & \exists \ C_0 \geq 0,  \forall x\in \partial \mathcal{O},
\forall x' \in \bar{\mathcal{O}} ,  \exists \  k \in {\bf
n}(x)\nonumber\\
&& \Longrightarrow  (x-x', k) + C_0|x-x'|^2\geq 0
,\\
& (ii) &  \forall  x\in  \partial \mathcal{O}, \mbox{if}   \exists C
\geq 0, \exists k \in \Re^d  , \forall  x' \in \bar{ \mathcal{O}}
, \nonumber\\
&& (x-x', k ) + C|x-x'|^2\geq 0 , \Longrightarrow  k =\theta {\bf
n}(x)
\end{eqnarray}
for some  $  \theta \geq 0$ , where $ \partial \mathcal{O} $ denotes
the boundary of $\mathcal{O}$, $\bar{\mathcal{O}} $ denotes the
closure of $\mathcal{O}$.
 We assume
that $ B_t $ is an  $ \Re^d $- valued $\mathcal{F}_t   $-  Brownian
motion on  a stochastic basis $(\Omega, {\mathcal F},\{{\mathcal
F}_t\}_{t\in [0,1]}, {\bf P})$ satisfying the usual assumptions.  We
consider the following stochastic differential equations on   domain
$ \mathcal{O} $ with reflecting boundary conditions:
\begin{eqnarray}
X_t(x)=x + \int^t_0 b(X_s)ds +\int^t_0 \sigma (X_s(x))\circ dB_s
-L^x_t, \ \forall \ t \in [0,1],
\end{eqnarray}
where $b: \Re^d\mapsto \Re^d $ and $ \sigma: \Re^d \mapsto \Re^d
\times \Re^d $ are continuous functions,  $\circ $ denotes the
Stratonovich integral. A pair  $ (X_t(x), L_t^x, t\in [0,1])$ of
continuous and $\mathcal{F}_t $- adapted processes  is called a
solution to equations (1.3) if  there exists a measurable set
$\tilde{ \Omega } $ with ${\bf P}( \tilde{
\Omega })=1$ such that for each $\omega \in \tilde{ \Omega } $\\
(i) for each $ x \in \bar{\mathcal{O}}$ the function $s \mapsto$ $
L_s^x$ with  values in $\Re^d $  has bounded
variation on any interval $ [0, T]$ and $ L_0^x=0 $. \\
(ii) for all $t\geq 0 $, $X_t(x)\in \bar{\mathcal{O}}$ and $
(X_t(x), L_t^x, t\in [0,1])$ satisfies Eq.(1.3).\\
(iii)
\begin{eqnarray}
|L^x|_t=\int^t_0 I_{(  X_s(x) \in  \partial \mathcal{O} ) }d|L^x|_s\
\mbox{ and }\  L_t^x=\int^t_0 {\bf \xi}(X_s(x))d|L^x|_s
\end{eqnarray}
with ${\bf \xi }(X_s(x)) \in {\bf n}(X_s(x)) $, where the $ |L^x|_t
$ denotes the total variation of $  L_\cdot^x$ on $[0, t]$.

Remark that  (iii) implies that the support of $d|L^x|$ is included
in $\{ s: X_s(x)\in
\partial \mathcal{O}\}$ and the force $ L^x $ keeps the process $ X $  be in
$\bar{\mathcal{O}}$.\\

This type of reflected stochastic differential equations has been
studied notably by Skorohod\cite{s1}, Tanaka\cite{s2}, Lions and
Sznitman\cite{s3},  and Saisho\cite{s4}, and also by Stroock and
Varadhan\cite{s5} who used a {\sl submartingale problem }
formulation, and other authors. Moreover, such reflected diffusions
can also be reduced to studying multivalued  stochastic differential
equations( see \cite{s6,s7,s8,s9} and  references therein).  It is
well-known (see \cite{s3}) that for any given initial value $x \in
\bar{\mathcal{O}}  $ the Eq.(1.3) has a unique solution
 provided that $\|\sigma(\cdot)\|$ and $|\widetilde{b}(\cdot)|$ are uniformly
bounded real-valued functions on $\Re^d$ and satisfy a uniform
Lipschitz condition: $\exists \ c>0  $ such that
\begin{eqnarray}
&&\|\sigma(y)-\sigma(z)\|\leq c|y-z|,\
|\widetilde{b}(y)-\widetilde{b}(z)|\leq c|y-z|
\end{eqnarray}
for any $y, z \in \Re^d$, where $ \widetilde{b}_i(x)=b_i(x)+
\frac{1}{2}\sum\limits^d_{k,j=1} \frac{\partial
\sigma_{ij}}{\partial x_k} (x)\sigma_{kj}(x)$,
$\|\sigma(y)\|:=\sqrt{\sum\limits^d_{i,j=1}\big\{ \sigma_{ij}(y)
\big\}^2 }$ and $|\widetilde{b}(y)|:=\sqrt{ \sum\limits^d_{i=1} \big
\{\widetilde{b}_i(y)\big\}^2}$. \\

{{\bf \sl The natural question aries: does there still exist a pair
$ (X_t, L_t, t\in [0,1])$ of stochastic processes to solve Eq.(1.3)
if the initial value is an arbitrary  random variable $Z$ which
belongs to $\bar{\mathcal{O}} $
 with probability one and may depend on the whole Brownian paths ?}}
 \vskip 0.3cm

 On one hand, the answer is
not immediately clear because one needs to deal with anticipating
stochastic integration. On the other hand, on a given financial
market, different agents generally have different levels of
information; besides the public information, some of them may
possess privileged information, which leads them to make
anticipations on some future realizations of functionals of the
price process, therefore, for a financial corporation, the studying
the problem of optimal dynamic risk control/dividends distribution
has to face the question(see \cite{s15, s16, s17, s18} and
references therein).  The main aim of this paper is to give an
affirmative answer to the  question above. Let us describe now more
precisely main results of this paper as follows.
\begin{Them}
Assume that $\mathcal{O} $ is a smooth bounded open set in $\Re^d$
and there exists a function $\phi\in {\mathcal{C}}_b^2(\Re^d ) $
such that
\begin{eqnarray}
\exists \alpha >0,\ \forall \in \partial \mathcal{O},\ \forall \zeta
\in {\bf n}(x),\ (\triangledown \phi(x), \zeta ) \leq -\alpha C_0,
\end{eqnarray}
the functions $\sigma $ and ${b} $ satisfy  that $\sigma $,
$\widetilde{b} $ and $ \triangledown \sigma     $ are bounded, and
  the following
\begin{eqnarray}
&&| \widetilde{b}(x)-\widetilde{b}(y)| + \|\sigma (x)-\sigma(y) \|+
\|(\triangledown \sigma \cdot \sigma)(x)- (\triangledown \sigma
\cdot \sigma )(y)\|\nonumber\\
&&\|(\triangledown \sigma  \cdot \triangledown \sigma \cdot
\sigma)(x)- (\triangledown \sigma \cdot \triangledown \sigma \cdot
\sigma )(y)\| + \| (\triangledown \sigma \cdot \widetilde{b})(x)-
(\triangledown \sigma \cdot \widetilde{b})(y)\|\nonumber\\
&&+\| (\sigma ^T\cdot\triangledown^2 \sigma \cdot \sigma)(x)-
(\sigma^T\cdot \triangledown^2 \sigma\cdot \sigma)(y)\|\leq k |x-y|
\end{eqnarray}
for some constant $k>0 $, where $C_0 $ is given by (1.1), $
\sigma^T$ denotes transpose of $\sigma$,  $\triangledown \sigma $
and  $ \triangledown^2 \sigma$ denote $\sigma$'s   derivatives of
first and second order with respect to spatial variable $x$,
respectively.
 Then for any random variable $Z $ with ${\bf P}\{ Z\in \bar{\mathcal{O}}
 \}=1 $
 the pair $ (X_t(Z), L_t^Z, t\in [0,1])$  solves the following
stochastic differential equation on domain $ \mathcal{O} $ with
reflecting boundary conditions:
\begin{eqnarray}
X_t(Z)=Z + \int^t_0 b(X_s(Z))ds +\int^t_0 \sigma (X_s(Z))\circ dB_s
-L^Z_t
\end{eqnarray}
with  $ X_t(Z)\in \bar{\mathcal{O}}$, and
  satisfies \\
(1)\ the function $s \mapsto$ $ L_s^Z$ with values in $\Re^d $ has
bounded variation on any interval $ [0, T]$ and $ L_0^Z=0 $. \\
(2)
\begin{eqnarray}
|L^Z|_t=\int^t_0 I_{(  X_s(Z) \in  \partial \mathcal{O} ) }d|L^Z|_s\
\mbox{ and }\  L_t^Z=\int^t_0 {\bf \xi}(X_s(Z))d|L^Z|_s
\end{eqnarray}
with ${\bf \xi }(X_s(Z)) \in {\bf n}(X_s(Z)) $, where  $ (X_t(x),
L^x_t, t\in [0,1])$ is the unique solution of Eq.(1.3),
 the stochastic integral in Eq.(1.8) is interpreted as
anticipating Stratonovich integral.
\end{Them}
Now we recall the definition of the anticipating Stratonovich
integral (see \cite{s10}). For any $t\in [0,1]$, let $\pi$ denote an
arbitrary partition of the interval $[0,t]$ of the form: $\pi =
\{0=t_0< t_1< \cdots < t_n=t \} $. Let $ \|\pi \|=\sup\limits_{0\leq
k \leq n-1}\{(t_{k+1}-t_k ) \}$ denote the norm of $\pi $. For an
$\Re^d \times \Re^d$-valued stochastic process $ f=\{ f_s, \ s\in
[0,1] \}$, we define its Riemann sums $S_\pi(f, t) $ by
\begin{eqnarray}
S_\pi(f, t)=\sum^{n-1}_{k=0} \frac{1}{t_{k+1}-t_k}\bigg (
\int^{t_{k+1}}_{t_k} f_s ds  \bigg )\cdot ( B_{t_{k+1}}-B_{t_k}).
\end{eqnarray}
We have the following
\begin{Def}
We say that a stochastic  process $ f=\{ f_s, \ s\in [0,1] \}$  is
Stratonovich integrable  with respect to $B$ if the family $S_\pi(f,
t) $ converges in probability as $\|\pi \|\rightarrow  0 $. In such
a case  we denote the limit  by $ \int^t_0 f_s\circ dB_s $.
\end{Def}
Let us now describe our approach. To prove Theorem 1.1, the natural
idea is to replace $x$ in (i), (ii) and (iii) of  Eq.(1.3) by the
initial random variable $Z$ and prove that the pair $( X_t(Z),
L^Z_t)$ solves  the Eq.(1.8). To achieve this, the key is to
establish the following
 substitution formula
\begin{eqnarray}
&& \int^t_0  \sigma (X_s(x))\circ dB_s\big |_{x=Z}=\int^t_0  \sigma
(X_s(Z))\circ  d B_s , \\
&& \int^t_0 f(  X_s(Z))d|L^Z|_s=0, \ \ L_t^Z=\int^t_0 {\bf
\xi}(X_s(Z))d|L^Z|_s
\end{eqnarray}
for all $t\in [0,1]$, where $f$ is a continuous function defined on
$\Re^d$ with compact support included in $\mathcal{O}$ and  ${\bf
\xi }(X_s(Z)) \in {\bf n}(X_s(Z)) $.\\

The novelty and difficulty of this paper are  anticipation,
reflection and shape of domain $ \mathcal{O}$. Since Lions and
Sznitman's result in \cite{s3} states  that the solution $ ( X_t(x),
L^x_t )$ is H\"{o}lder continuous of order being less than
$\frac{1}{2}$ with respect to the initial value $x$, the regularity
 is not good enough to satisfy the required hypothesis of substitution
 formula in the literature (see \cite{s10}), it seems that we can
  not apply the existing
substitution formula to prove (1.11). Moreover, because reflecting
boundary conditions and shape of domain $ \mathcal{O}$, it  is  also
impossible  to prove (1.11) by using It\^{o}-Ventzell formula used
by cone and Pardoux\cite{s19}, Kohatsu-Higa and Le\'{o}n\cite{s20}.
 Instead, we prove (1.11) by showing  the uniform
convergence (w.r.t.x) of the corresponding Riemann Sums $S_\pi(
\sigma (X_\cdot(x),t)$. The Garsia, Rodemich and Rumsey's Lemma and
moments estimates for one-point and two-point motions will play an
important role. To prove (1.12) we need only to show that  the
functionals $ F(t,x):=\int^t_0 f( X_s(x))d|L^x|_s$, $L^x_t$ and $
G(t,x):=\int^t_0 {\bf \xi}(X_s(x))d|L^x|_s$ are continuous in
$(t,x)$, doing this depends on the shape function $\phi  $ of domain
$ \mathcal{O}$ in (1.6).
\begin{remark}
It seems  that this new approach can be used to study
 perturbed stochastic Skorohod equations with anticipating initial
 processes because the solutions of these SDEs are H\"{o}lder continuous
 with order being
less than $\frac{1}{2}$ and not differentiable w.r.t. initial value
$x$ (see \cite{s11,s12}for adapted case). We shall study it in
forthcoming paper.
\end{remark}
 This paper is organized as follows. Firstly we study the regularity
 of the solution $ ( X_t(x), L^x_t )$. Secondly we devote to showing
continuity of functionals $ F(t,x):=\int^t_0 f( X_s(x))d|L^x|_s$,
$L^x_t$ and $ G(t,x):=\int^t_0 {\bf \xi}(X_s(x))d|L^x|_s$ . In
Section 4 we study  moments estimates for one-point and two-point
motions. In Section 5 we prove the uniform convergence (w.r.t.x) of
the Riemann Sums $S_\pi( \sigma (X_\cdot(x),t)$. Finally we prove
Theorem 1.1 in Section 6.\\

Throughout this paper we make the following convention: the letter
$c$ or $c(p_1, p_2, p_3, \cdot, \cdot, \cdot )$ depending only on $
p_1, p_2, p_3, \cdot, \cdot, \cdot $ will denote an unimportant
positive constant, whose values may change from one line to another
one.
\setcounter{equation}{0}
\section{\bf {\small Regularity of the solution $ ( X_t(x), L^x_t ) $ of Eq.(1.3)}}
The main aim of this section is  to study regularity of the solution
$ ( X_t(x), L^x_t )$ $ w.r.t. (t, x)$ via the shape function $\phi $
of domain $ \mathcal{O}$ in (1.6).
\begin{Prop}
Assume that the smooth bounded open $  \mathcal{O} $, the
coefficients $ \sigma $ and $b$ satisfy  the same conditions as in
Theorem 1.1.  $ ( X_t(x), L^x_t ) $ is a solution of Eq.(1.3). Then
there is a constant $c$ such that
\begin{eqnarray}
&&{\bf E}\big \{\sup_{0\leq t\leq 1}|X_t(x)-X_t(y)|^{p}  \big \}
\leq c|x-y|^p,\\
&&{\bf E}\big \{\sup_{0\leq t\leq 1}|L_t^x-L_t^y|^{p}  \big \} \leq
c|x-y|^p
\end{eqnarray}
for any $x, y \in  \bar{\mathcal{O} } $ and $p\geq 1$.
\end{Prop}
{\bf Proof.}\ By H\"{o}lder inequality, we need only to prove
Proposition 2.1 for $p\geq 4$. Let $ \widetilde{b}_i(x)=b_i(x)+
\frac{1}{2}\sum\limits^d_{k,j=1} \frac{\partial
\sigma_{ij}}{\partial x_k} (x)\sigma_{kj}(x)$, that is, $
\widetilde{b}(x)= b(x) + \frac{1}{2} (\triangledown \sigma \cdot
\sigma )(x) $   for any $x\in \Re^d $. We write solution $ ( X_t(x),
L^x_t )$ of Eq.(1.3) in It\^{o}'s form as follows: for $x\in
\bar{\mathcal{O} }$
\begin{eqnarray}
&&X_t(x)=x + \int^t_0 \widetilde{b}(X_s)ds +\int^t_0 \sigma (X_s(x))
dB_s -L^x_t, \\
&&|L^x|_t=\int^t_0 I_{(  X_s(x) \in  \partial \mathcal{O} )
}d|L^x|_s,\\
&&L_t^x=\int^t_0 {\bf \xi}(X_s(x))d|L^x|_s \ \mbox{ with}\  {\bf \xi
}(X_s(x)) \in {\bf n}(X_s(x))
\end{eqnarray}
and $ ( X_t(y), L^y_t )$ also satisfy the same equations above for
$y\in \bar{\mathcal{O} }$.

Applying It\^{o}'s formula to function $\phi \in \mathcal{C}_b^2
(\Re^d)$ satisfying (1.6) and stochastic process $X_t(x)$,  we have
\begin{eqnarray}
\phi(X_t(x))&=& \phi(x) + \int^t_0 ( \triangledown \phi^T \sigma )(
X_s(x))dB_s + \int^t_0 ( \triangledown \phi^T \widetilde{b} )(
X_s(x))ds\nonumber\\
&&- \int^t_0 ( \triangledown \phi^T \xi )( X_s(x))d|L^x|_s
\nonumber\\
&& + \frac{1}{2}\int^t_0 {\bf tr} \big \{\big (\triangledown^2\phi
\sigma \sigma^T \big )(X_s(x)) \big\}ds,
\end{eqnarray}
where ${\bf tr(A) }$ denote the trace of $ A$.
Similarly, we have same expression for $\phi(X_t(y))$.\\
Define $ f(x):= |x|^p $, $x=(x_1, x_2, \cdots, x_d)^T \in \Re^d$.
Then
\begin{eqnarray}
\triangledown f(x)=p|x|^{p-2}x, \ \triangledown^2
f(x)=p|x|^{p-2}I_{d\times d} + p(p-2)|x|^{p-4}xx^T.
\end{eqnarray}
Let $m_t= X_t(x)-X_t(y)$, $D_t=\phi(X_t(x)) + \phi(X_t(y))$ and
$N_t=\exp\{-\frac{p}{\alpha}D_t \}$.
 By It\^{o}'s formula and (2.5),
\begin{eqnarray}
df(m_t)&=& p|m_t|^{p-2}m_t^T(\widetilde{b}(X_t(x))
-\widetilde{b}(X_t(y))) dt\nonumber \\
&+& p|m_t|^{p-2}m_t^T(\sigma(X_t(x))
-\sigma(X_t(y))) dB_t\nonumber \\
&-& p|m_t|^{p-2}m_t^T{\bf \xi}(X_t(x))d|L^x|_t\nonumber\\
&+& p|m_t|^{p-2}m_t^T{\bf \xi}(X_t(y))d|L^y|_t\nonumber\\
&+&\frac{1}{2}{\bf tr} \big\{\triangledown^2f(m_t) (\sigma(X_t(x))
-\sigma(X_t(y))) (\sigma(X_t(x)) -\sigma(X_t(y)))^T \}dt,\nonumber\\
\end{eqnarray}
\begin{eqnarray}
dN_t&=& -\frac{p}{\alpha}N_t\big [(\triangledown \phi^T \sigma)
(X_t(x))+ ( \triangledown \phi^T \sigma )(X_t(y)) \big
]dB_t\nonumber\\
&&-\frac{p}{\alpha}N_t\big [(\triangledown \phi^T \widetilde{b})
(X_t(x))+ ( \triangledown \phi^T \widetilde{b} )(X_t(y)) \big
]dt\nonumber\\
&&+\frac{p}{\alpha}N_t(\triangledown \phi^T {\bf\xi })
(X_t(x))d|L^x|_t\nonumber\\
&&+\frac{p}{\alpha}N_t(\triangledown \phi^T {\bf\xi })
(X_t(y))d|L^y|_t\nonumber\\
&&-\frac{p}{2\alpha}N_t {\bf tr}\big\{ \big (\triangledown^2\phi
\sigma \sigma^T \big )(X_t(x)) +\big (\triangledown^2\phi \sigma
\sigma^T
\big )(X_t(y)) \big\}dt\nonumber\\
&&+\frac{p^2}{2\alpha^2} N_t {\bf tr} \bigg \{\big [(\triangledown
\phi^T \sigma) (X_t(x))+ ( \triangledown \phi^T \sigma )(X_t(y))
\big
]^T\nonumber\\
&&\times\big [(\triangledown \phi^T \sigma) (X_t(x))+ (
\triangledown \phi^T \sigma )(X_t(y)) \big ] \bigg \}dt
\end{eqnarray}
and the stochastic contraction $ df(m_t)\cdot dN_t $ is given by
\begin{eqnarray}
df(m_t)\cdot dN_t &=&-\frac{p^2}{\alpha} N_t  |m_t|^{p-2} {\bf tr}
\bigg\{
\big(m_t^T(\sigma(X_t(x)) -\sigma(X_t(y))) \big )^T\nonumber\\
&& \times\big [(\triangledown \phi^T \sigma) (X_t(x))+ (
\triangledown \phi^T \sigma )(X_t(y)) \big ]^T \bigg\}dt.
\end{eqnarray}
Therefore, by It\^{o}'s formula again and (2.8)-(2.10),
\begin{eqnarray}
N_tf(m_t)&=& \exp\{-\frac{p}{\alpha}[\phi(x)+\phi(y) ] \}|x-y|^p
+\int^t_0 N_sdf(m_s)\nonumber\\
 &+& \int^t_0 f(m_s)dN_s +\int^t_0
df(m_s)\cdot dN_s\nonumber\\
&=&\exp\{-\frac{p}{\alpha}[\phi(x)+\phi(y) ] \}|x-y|^p\nonumber\\
&+& p\int^t_0N_s|m_s|^{p-2}m_s^T(\widetilde{b}(X_s(x))
-\widetilde{b}(X_s(y))) ds\nonumber \\
&+& p\int^t_0N_s|m_s|^{p-2}m_s^T(\sigma(X_s(x))
-\sigma(X_s(y))) dB_s\nonumber \\
&-& p\int^t_0N_s|m_s|^{p-2}\big ( m_s, {\bf \xi}(X_t(x))\big )d|L^x|_s\nonumber\\
&+& p\int^t_0N_s|m_s|^{p-2}\big ( m_s, {\bf \xi}(X_t(y))\big )d|L^y|_s\nonumber\\
&+&\frac{1}{2}\int^t_0N_s {\bf tr}\big\{\triangledown^2f(m_s)
(\sigma(X_s(x))-\sigma(X_s(y))) (\sigma(X_s(x))\nonumber\\
& -&\sigma(X_s(y)))^T \}ds\nonumber\\
&-&\frac{p}{\alpha}\int^t_0 N_sf(m_s)\big [(\triangledown \phi^T
\sigma) (X_s(x))+ ( \triangledown \phi^T \sigma )(X_s(y)) \big
]dB_s\nonumber\\
&-&\frac{p}{\alpha}\int^t_0 N_sf(m_s)\big [(\triangledown \phi^T
\widetilde{b}) (X_s(x))+ ( \triangledown \phi^T \widetilde{b}
)(X_s(y)) \big]ds\nonumber\\
&+&\frac{p}{\alpha}\int^t_0 N_s|m_s|^{p-2}|m_s|^2 \big(\triangledown
\phi(X_s(x)),  {\bf\xi }(X_s(x))\big )d|L^x|_s
\nonumber\\
&+& \frac{p}{\alpha}\int^t_0 N_s|m_s|^{p-2}|m_s|^2
\big(\triangledown \phi(X_s(y)),  {\bf\xi }(X_s(y))\big )d|L^y|_s
\nonumber\\
&-&\frac{p}{2\alpha}\int^t_0 N_sf(m_s){\bf  tr}\big\{ \big
(\triangledown^2\phi \sigma \sigma^T \big )(X_s(x)) +\big
(\triangledown^2\phi \sigma \sigma^T
\big )(X_s(y)) \big\}ds\nonumber\\
&+&\frac{p^2}{2\alpha^2} \int^t_0 N_sf(m_s){\bf  tr}\bigg \{\big
[(\triangledown \phi^T \sigma) (X_s(x))+ ( \triangledown \phi^T
\sigma )(X_s(y)) \big
]^T\nonumber\\
&&\times\big [(\triangledown \phi^T \sigma) (X_s(x))+ (
\triangledown \phi^T \sigma )(X_s(y)) \big ] \bigg \}ds\nonumber\\
&-&\frac{p^2}{\alpha}\int^t_0 N_s  |m_s|^{p-2} {\bf tr } \bigg\{
\big(m_s^T(\sigma(X_s(x)) -\sigma(X_s(y))) \big )^T\nonumber\\
&& \times\big [(\triangledown \phi^T \sigma) (X_s(x))+ (
\triangledown \phi^T \sigma )(X_s(y)) \big ]^T \bigg\}ds\nonumber\\
&:=&\sum_{i=1}^{13} a_i(t).
\end{eqnarray}
By condition (1.6),
\begin{eqnarray*}
&&\frac{1}{\alpha}|m_s|^2 \big(\triangledown \phi(X_s(x)), {\bf\xi
}(X_s(x))\big )-\big ( m_s, {\bf \xi}(X_t(x))\big )\leq 0,\ \
d|L^x|_s\ a.s.\\
&&\frac{1}{\alpha}|m_s|^2 \big(\triangledown \phi(X_s(y)), {\bf\xi
}(X_s(y))\big )+\big ( m_s, {\bf \xi}(X_t(y))\big )\leq 0,\ \
d|L^y|_s\ a.s.\\
\end{eqnarray*}
Hence
\begin{eqnarray}
a_4(t)+a_9(t)\leq 0, \quad  a_5(t) +a_{10}(t)\leq 0.
\end{eqnarray}
Using $ \phi \in \mathcal{C}_b^2(\Re^d)$ and (2.11)-(2.12),
\begin{eqnarray}
(m^*(t))^{2p}\leq c [
\sum_{i=1}^3(a^*_i(t))^2+\sum_{i=6}^8(a^*_i(t))^2
+\sum_{i=11}^{13}(a^*_i(t))^2],
\end{eqnarray}
where $a^*_i(t)=\sup_{s\in[0,t]}\{|a_i(t)|\}$,
$m^*(t)=\sup_{s\in[0,t]}\{|m(t)|\}$.

Since $ \phi $ is bounded, by Burkh\"{o}lder inequality(see
\cite{s13}) and (1.7), we have
\begin{eqnarray}
{\bf E} \{ (a_3^*(t))^2\} &\leq & c{\bf E}\bigg \{ \int^t_0
N_s^2|m_s|^{2p-4}{\bf tr} \big\{ [ m_s^T(\sigma(X_s(x))
-\sigma(X_s(y)))
]^T\nonumber\\
&&\times [m_s^T(\sigma(X_s(x)) -\sigma(X_s(y))) ] \big\}
ds\bigg\}\nonumber\\
&\leq & c{\bf E}\bigg \{ \int^t_0 |m_s|^{2p-2} \|(\sigma(X_s(x))
-\sigma(X_s(y)))\|^2ds\bigg\}\nonumber\\
&\leq & c \int^t_0{\bf E} \{ |m^*(s)|^{2p}\} ds.
\end{eqnarray}
Similarly, since $ \phi $  and  $\triangledown \phi \ \sigma    $
are bounded on $  \bar{\mathcal{O}}$, we also have
\begin{eqnarray}
{\bf E} \{ (a_7^*(t))^2\} &\leq &  c \int^t_0{\bf E} \{
|m^*(s)|^{2p}\} ds.
\end{eqnarray}
Using $ \phi $  and  $\triangledown \phi \ \sigma $ are bounded on $
\bar{\mathcal{O}}$, the condition (1.7) and H\"{o}lder inequality,
\begin{eqnarray}
{\bf E} \{ (a_{13}^*(t))^2\} &\leq & c{\bf E}\bigg \{ \int^t_0
N_s^2|m_s|^{2p-4}\big | {\bf tr} \bigg\{
\big(m_s^T(\sigma(X_s(x)) -\sigma(X_s(y))) \big )^T\nonumber\\
&& \times\big [(\triangledown \phi^T \sigma) (X_s(x))+ (
\triangledown \phi^T \sigma )(X_s(y)) \big ]^T \bigg\}\big |^2ds\nonumber\\
&\leq & c{\bf E}\bigg \{ \int^t_0 |m_s|^{2p-2} \|(\sigma(X_s(x))
-\sigma(X_s(y)))\|^2ds\bigg\}\nonumber\\
&\leq & c \int^t_0{\bf E} \{ |m^*(s)|^{2p}\} ds.
\end{eqnarray}
By the same way as in (2.16)
\begin{eqnarray}
{\bf E} \{ (a_i^*(t))^2\} &\leq &  c \int^t_0{\bf E} \{
|m^*(s)|^{2p}\} ds
\end{eqnarray}
 for $i=2,6,8,11,12   $.\\
Putting the above inequalities (2.13)-(2.17)together implies  that
\begin{eqnarray}
{\bf E} \{ (m^*(t))^{2p}  \} &\leq & c|x-y|^{2p}+ c \int^t_0{\bf E}
\{ (m^*(s))^{2p}\} ds,
\end{eqnarray}
By Gronwall inequality,
\begin{eqnarray}
{\bf E} \{ \sup_{t\in [0,1]}  |X_t(x)-X_t(y)|^{2p}  \} &\leq &
c|x-y|^{2p}.
\end{eqnarray}
So proof of (2.1) has been done by  H\"{o}lder inequality. Using
\begin{eqnarray*}
&&|L^x_t-L^y_t|\leq |x-y| +|X_t(x)-X_t(y)| +
|\int^t_0(\widetilde{b}(X_s(x))
-\widetilde{b}(X_s(y))) ds|\\
&&+ |\int^t_0(\sigma(X_s(x)) -\sigma(X_s(y))) dB_s|,
\end{eqnarray*}
(2.2) is a direct consequence of (2.1). Thus we complete  proof of
Proposition 2.1. \quad $\Box$
\begin{Prop}
Assume that the smooth bounded open $  \mathcal{O} $, the
coefficients $ \sigma $ and $b$ satisfy  the same conditions as in
Theorem 1.1.  $ ( X_t(x), L^x_t ) $ is a solution of Eq.(1.3). Then
there is a constant $c$, which is independent of $x$, such that for
any $p\geq 2$
\begin{eqnarray}
&&\sup_{x\in \bar{\mathcal{O}}    }{\bf E}\big
\{|X_t(x)-X_s(x)|^{2p} \big \}
\leq c|t-s|^{\frac{p}{2}},\\
&&\sup_{x\in \bar{\mathcal{O}}    }{\bf E}\big \{|L_t^x-L_s^x|^{2p}
\big \}\leq c|t-s|^{\frac{p}{2}}.
\end{eqnarray}
\end{Prop}
{\bf Proof.}\ By H\"{o}lder inequality, we need only to prove
Proposition 2.2 for $p\geq 4$. For $t\geq s\geq 0 $, similar to that
of Proposition 2.1,  we define $m_t$,$D_t$, $N_t$ and $f$ here  by
\begin{eqnarray*}
m_t(x)&=&X_t(x)-X_s(x)\\
&=&\int^t_s \widetilde{b}(X_u(x))du +\int^t_s \sigma (X_u(x))
dB_u - \int^t_s {\bf \xi}(X_u(x))d|L^x|_u\nonumber\\
&& \mbox{ with}\  {\bf \xi
}(X_s(x)) \in {\bf n}(X_s(x)), \\
D_t&=&\phi(X_t(x)),\\
 N_t&=&\exp\{-\frac{p}{\alpha}D_t \},\\
G_t&= &N_t^{-1},\\
f(x)&=& |x|^2 , x=(x_1, x_2, \cdots, x_d)^T \in \Re^d.
\end{eqnarray*}
By the same way as in (2.11),
\begin{eqnarray*}
|X_t(x)-X_s(x)|^2&=& 2G_t\int^t_sN_u m_u^T\widetilde{b}(X_u(x))du\nonumber \\
&+& 2G_t\int^t_sN_um_u^T\sigma(X_u(x))dB_u\nonumber \\
&-& 2G_t\int^t_sN_u\big ( m_u, {\bf \xi}(X_u(x))\big )d|L^x|_u\nonumber\\
&+&G_t\int^t_sN_u{\bf tr}\big\{ (\sigma \sigma^T)(X_u(x)) \}du\nonumber\\
&-&\frac{2}{\alpha}G_t\int^t_s N_uf(m_u)\big [(\triangledown \phi^T
\sigma) (X_u(x))\big ]dB_u\nonumber\\
&-&\frac{2}{\alpha}G_t\int^t_s N_uf(m_u)(\triangledown \phi^T
\widetilde{b}) (X_u(x))du\nonumber\\
&+&\frac{2}{\alpha}G_t\int^t_s N_u|m_u|^2 \big(\triangledown
\phi(X_u(x)),  {\bf\xi }(X_u(x))\big )d|L^x|_u
\nonumber\\
&-&\frac{1}{\alpha}G_t\int^t_s N_uf(m_u){\bf  tr}\big\{ \big
(\triangledown^2\phi \sigma \sigma^T \big )(X_u(x)) \big\}du\nonumber\\
&+&\frac{2}{\alpha^2}G_t \int^t_s N_uf(m_u) {\bf tr}\bigg \{\big
[(\triangledown \phi^T \sigma) (X_u(x)) \big]^T  \big
[(\triangledown
\phi^T \sigma) (X_u(x)) \big ] \bigg \}du\nonumber\\
&-&\frac{4}{\alpha}G_t\int^t_s N_u {\bf tr} \bigg\{
\big(m_u^T(\sigma(X_u(x))) \big )^T\big ((\triangledown \phi^T
\sigma) (X_u(x)) \big )^T \bigg\}du\nonumber\\
&:=&\sum_{i=1}^{10} d_i(t).
\end{eqnarray*}
By condition (1.6),
\begin{eqnarray}
d_3(t)+ d_7(t)\leq 0.
\end{eqnarray}
Therefore
\begin{eqnarray}
{\bf E}\{|X_t(x)-X_s(x)|^{2p}\}&\leq & c(p)\sum^2_{i=1}{\bf E}\{|
d_i(t) |^p \} + c(p)\sum^6_{i=4}{\bf E}\{| d_i(t) |^p \} \nonumber\\
&&+ c(p)\sum^{10}_{i=7}{\bf E}\{| d_i(t) |^p \}.
\end{eqnarray}
Since $ \sigma $, $ N_t$ and $G_t$ are uniformly bounded, by
Burkh\"{o}lder (see \cite{s13}) and H\"{o}lder inequalities and
Young's inequality: for any real positive $x,y, \eta, p, q$ with $
p^{-1} + q^{-1}=1$ there exists $c< +\infty $ such that $ xy\leq
\eta x^p +cy^q $, we have
\begin{eqnarray}
{\bf E}\{ |d_2(t)|^p\}&\leq & c{\bf E}\big \{ | \int^t_s N_u m_u^T
\sigma(X_u(x))dB_u|^p\big \}\nonumber\\
&\leq & c {\bf E}\big \{ \int^t_s N_u^2 {\bf
tr}\{\sigma^T(X_u(x))m_u
m_u^T \sigma (X_u(x))\}du\big \}^{\frac{p}{2}}\nonumber\\
&\leq & c {\bf E}\big \{ \int^t_s |m_u|^2du\big \}^{\frac{p}{2}}\nonumber\\
&\leq & c|t-s|^{\frac{p}{2}}+c \int^t_s{\bf E}\{ |m_u|^{2p}\}du.
\end{eqnarray}
Similarly,
\begin{eqnarray}
{\bf E}\{ |d_5(t)|^p\}&\leq & c \int^t_s{\bf E}\{ |m_u|^{2p}\}du.
\end{eqnarray}
Since $ \sigma $, $\triangledown \phi  \sigma  $, $ N_t$ and $G_t$
are uniformly bounded, by H\"{o}lder inequalities and Young's
inequality, we have
\begin{eqnarray}
{\bf E}\{ |d_{10}(t)|^p\}&\leq &
 c{\bf E}\bigg \{ |\int^t_s N_u {\bf tr} \bigg\{
\big(m_u^T(\sigma(X_u(x))) \big )^T\big ((\triangledown \phi^T
\sigma) (X_u(x)) \big )^T \bigg\}du |^p\bigg \}\nonumber\\
&\leq & c  {\bf E}\bigg \{ \int^t_s |m_u|du\bigg \}^{p}\nonumber\\
&\leq & c|t-s|^{p}+c \int^t_s{\bf E}\{ |m_u|^{2p}\}du.
\end{eqnarray}
Similarly,
\begin{eqnarray}
&&{\bf E}\{ |d_{1}(t)|^p\}\leq  c|t-s|^{p}+c \int^t_s{\bf E}\{
|m_u|^{2p}\}du,\\
&&{\bf E}\{ |d_{4}(t)|^p\}\leq c|t-s|^p,\\
&&{\bf E}\{ |d_{i}(t)|^p\}\leq  c \int^t_s{\bf E}\{ |m_u|^{2p}\}du,
\quad i=6,8,9.
\end{eqnarray}
Putting the above inequalities (2.23)-(2.29) together, we obtain
\begin{eqnarray*}
{\bf E}\{|X_t(x)-X_s(x)|^{2p}\}\leq  c|t-s|^{\frac{p}{2}}+c
\int^t_s{\bf E}\{ | X_u(x)-X_u(x)     |^{2p}\}du.\nonumber\\
\end{eqnarray*}
The Gronwall-Bellman inequality(see \cite{s14} for Theorem 1.3.1)
implies that
$$ {\bf E}\{|X_t(x)-X_s(x)|^{2p}\}\leq  c|t-s|^{\frac{p}{2}}. $$
Therefore the  proof of (2.20) has been done. Using
\begin{eqnarray*}
&&|L^x_t-L^x_s|\leq |X_t(x)-X_s(x)| + |\int^t_s
\widetilde{b}(X_u(x))du | + |\int^t_s \sigma (X_u(x)) dB_u |,
\end{eqnarray*}
(2.21) is a direct consequence of (2.20). Thus we complete  proof of
Proposition 2.2. \quad $\Box$

Since the domain $ \mathcal{O}$ is bounded, the following follows
immediately from Proposition 2.1 and H\"{o}lder's inequality.
\begin{Prop}
Assume that the smooth bounded open $  \mathcal{O} $, the
coefficients $ \sigma $ and $b$ satisfy  the same conditions as in
Theorem 1.1.  $ ( X_t(x), L^x_t ) $ is a solution of Eq.(1.3). Then
there is a constant $c$, which is independent of $x$,  such that
\begin{eqnarray}
&&{\bf E}\big \{\sup_{0\leq t\leq 1}|X_t(x)|^{p}  \big \} \leq c(1+|x|)^p,\\
&&{\bf E}\big \{\sup_{0\leq t\leq 1}|L^x_t|^{p}  \big \} \leq
c(1+|x|)^p
\end{eqnarray}
for any $x\in  \bar{\mathcal{O} } $ and $p\geq 1$.
\end{Prop}
\vskip 0.3cm
\setcounter{equation}{0}
\section{ \bf {\small Continuity of functionals of local times} }
\begin{Prop}
Assume that the smooth bounded open $  \mathcal{O} $, the
coefficients $ \sigma $ and $b$ satisfy  the same conditions as in
Theorem 1.1.  $ ( X_t(x), L^x_t ) $ is a solution of Eq.(1.3). Then
the functions  $ X_t(x)$, $L^x_t$, $F(t,x)$ and $ G(t,x)$ are
jointly continuous in $(t, x)$ on $[0,1]\times \bar{\mathcal{O}}$,
where  $ F(t,x):=\int^t_0 f( X_s(x))d|L^x|_s$ and $ G(t,x):=\int^t_0
{\bf \xi}(X_s(x))d|L^x|_s$, $f$ is a continuous function defined on
$\Re^d$ with compact support included in $\mathcal{O}$ and  ${\bf
\xi }(X_s(x)) \in {\bf n}(X_s(x)) $.
\end{Prop}
{\bf Proof.} By Kolmogorov's continuity criterion of random
fields(see \cite{s21} for Theorem 1.4.1 ), Proposition 2.1 and 2.2,
the functions $ X_t(x)$ and $L^x_t$ are H\"{o}lder continuous in
$(t,x)$. Since proof of continuity of  $G(t,x)$  w.r.t.$(t,x)$ is
similar to that of $G(t,x)$, we need only deal with the proof of
$F(t,x)$. Remarking that
\begin{eqnarray}
|F(t,x)-F(s,x)|\leq \sup_{y\in
\bar{\mathcal{O}}}\{|f(y)|\}|L^x_t-L^x_s|,
\end{eqnarray}
the function $F(t,x)$ is continuous in $t$ uniformly with respect to
$x$ in  compact set $\bar{ \mathcal{O}}$ by Proposition 2.2 and
Kolmogorov's continuity criterion( see Theorem 1.4.1 in \cite{s21}
). Thus, it suffices to show the continuity of $ F(t,x)$ w.r.t.$x$
for any fixed $t$. Let $x_n,x\in \bar{ \mathcal{O}} $ with  $ x_n
\longrightarrow x $ as $ n\longrightarrow +\infty $. By Propositions
2.1-2.2, and  $ X_t(x)$ and $L^x_t$ are H\"{o}lder continuous in
$(t,x)\in [0,1]\times \bar{\mathcal{O}}$,  we have
\begin{eqnarray}
L^{x_n}_t\longrightarrow L^x_t, \quad  X_t(x_n)\longrightarrow
X_t(x),
\end{eqnarray}
uniformly in $t$, as $n\longrightarrow +\infty $. Therefore, there
exist  constants  $C_1,$ $  C\geq 1$ such that for all  $n\geq 1$
\begin{eqnarray}
 |L^{x_n}_t| \leq C +| L^x_t|\leq C + | L^x|_1\leq C + C_1
 \end{eqnarray}
due to bound of total variation of $L^x_\cdot$ on $[0,1]$.  Since
the function $f(x)$ is bounded and continuous, by (3.2) and (3.3),
\begin{eqnarray}
| \int^t_0 [ f(X_s(x_n))-f(X_s(x))]dL^{x_n}_s|\longrightarrow 0
\end{eqnarray}
as $n\longrightarrow +\infty$. Because $L^{x_n}_t $ and $ L^{x}_t$
are  continuous processes with bounded variation, by (3.2), the
sequence of finite  sign measures $ dL^{x_n}_t $ on $[0,1]$
converges weakly to the finite sign measure $ dL^{x}_t $ on $[0,1]$.
Therefore, for bounded continuous function $ f(X_s(x))$ on $[0,1]$,
we have
\begin{eqnarray}
\lim_{n\rightarrow \infty }\int^t_0 f( X_s(x))dL^{x_n}_s=\int^t_0 f(
X_s(x))dL^{x}_s.
\end{eqnarray}
The proof of Proposition 3.1 follows from (3.4) and (3.5).\quad
$\Box $\vskip 0.3cm
As a direct consequence of Proposition 3.1, we
have the following.
\begin{Prop}
Assume that the smooth bounded open $  \mathcal{O} $, the
coefficients $ \sigma $ and $b$ satisfy  the same conditions as in
Theorem 1.1.  $ ( X_t(x), L^x_t ) $ is a solution of Eq.(1.3). Then
there exists a set $ \widetilde{\Omega }\in \mathcal{F} $ with ${\bf
P} (\widetilde{\Omega })=1$ such that for each $\omega \in
\widetilde{\Omega }$
\begin{eqnarray}
&& \int^t_0  b(X_s(x))ds\big |_{x=Z}=\int^t_0 b(X_s(Z))ds , \\
&& \int^t_0 f(  X_s(Z))d|L^Z|_s=0, \ \ L_t^Z=\int^t_0 {\bf
\xi}(X_s(Z))d|L^Z|_s
\end{eqnarray}
for all $t\in [0,1]$, where $f$ is a continuous function defined on
$\Re^d$ with compact support included in $\mathcal{O}$ and  ${\bf
\xi }(X_s(Z)) \in {\bf n}(X_s(Z)) $.
\end{Prop}
\setcounter{equation}{0}
\section{{\bf {\small  Moments estimates for
one-point and two-point motions}}}
 For any $R >0 $ and  $ x \in [-R,
R]^d\cap \bar{\mathcal{O}}$, let $ ( X_t(x), L^x_t ) $
 be a solution of Eq.(1.3). We define  $ S_\pi (t,x) $ and $I(t,x) $
by\begin{eqnarray*}
S_\pi (t,x)&:=&S_\pi(\sigma (X_{\cdot} (x)), t),\\
I(t, x)&:=& \int^t_0\sigma (X_s (x))\circ dB_s\\
       & =&\int^t_0\sigma (X_s (x))dB_s +\frac{1}{2}
        \int^t_0 ( \triangledown \sigma   \cdot  \sigma )(X_s (x)) ds.
\end{eqnarray*}
Write
\begin{eqnarray}
&&S_\pi(t,x)\nonumber\\
&=&\sum^{n-1}_{k=0} \frac{1}{t_{k+1}-t_k}\bigg (
\int^{t_{k+1}}_{t_k} \sigma ( X_s(x) ) ds  \bigg )(
B_{t_{k+1}}-B_{t_k})\nonumber\\
&=& \sum^{n-1}_{k=0}\sigma ( X_{t_k}(x) )(
B_{t_{k+1}}-B_{t_k})\nonumber\\
&+&\sum^{n-1}_{k=0} \frac{1}{t_{k+1}-t_k}\bigg (
\int^{t_{k+1}}_{t_k} (\sigma ( X_s(x) )-\sigma ( X_{t_k}(x) )) ds
\bigg )( B_{t_{k+1}}-B_{t_k}).\nonumber\\
\end{eqnarray}
By Ito's formula and (1.4), for $s\geq t_k$,
\begin{eqnarray*}
\sigma_{ij} ( X_s(x) )-\sigma_{ij} ( X_{t_k}(x) ) &=& \int_{t_k}^s
\big (\triangledown
\sigma_{ij}\cdot \sigma \big )( X_u(x) ) dB_u\nonumber\\
&+& \int_{t_k}^s \big (\triangledown \sigma_{ij}\cdot
{b}\big ) ( X_u(x) )du\nonumber\\
& -& \int_{t_k}^s \big (\triangledown \sigma_{ij}\cdot {\bf \xi}\big
) ( X_u(x) )d|L^x|_u\nonumber\\
& + & \frac{1}{2}\int_{t_k}^s \bigg \{ \big (\triangledown
\sigma_{ij} \cdot (\triangledown \sigma \cdot\sigma) \big )( X_u(x)
)\bigg\}du\nonumber\\
 & + & \frac{1}{2}\int_{t_k}^s {\bf tr}\bigg \{ \big
(\triangledown^2 \sigma_{ij}\cdot
 \sigma \sigma^T \big )( X_u(x) )\bigg\}du.
\end{eqnarray*}
So we informally write $\sigma ( X_s(x) )-\sigma ( X_{t_k}(x) ) $ as
follows:
\begin{eqnarray}
\sigma ( X_s(x) )-\sigma ( X_{t_k}(x) ) &=& \int_{t_k}^s \big
(\triangledown
\sigma\cdot \sigma \big )( X_u(x) ) dB_u\nonumber\\
&+& \int_{t_k}^s \big (\triangledown \sigma\cdot
{b}\big ) ( X_u(x) )du\nonumber\\
& -& \int_{t_k}^s \big (\triangledown \sigma\cdot {\bf \xi}\big
) ( X_u(x) )d|L^x|_u\nonumber\\
& + & \frac{1}{2}\int_{t_k}^s \bigg \{ \big (\triangledown \sigma
\cdot (\triangledown \sigma \cdot\sigma) \big )( X_u(x)
)\bigg\}du\nonumber\\
 & + & \frac{1}{2}\int_{t_k}^s {\bf tr}\bigg \{ \big
(\triangledown^2 \sigma\cdot
 \sigma \sigma^T \big )( X_u(x) )\bigg\}du.\nonumber\\
\end{eqnarray}
Thus we can write $S_\pi(t,x)-I(t,x )$ as follows:
\begin{eqnarray}
S_\pi(t,x)-I(t,x )= \sum^6_{i=1}A_{i\pi},
\end{eqnarray}
where
\begin{eqnarray*}
A_{1\pi}(x)&:=&\sum^{n-1}_{i=0} \sigma
(X_{t_i}( x))( B_{t_{i+1}}-B_{t_i})-\int^t_0 \sigma(X_s(x))dB_s,\\
A_{2\pi}(x)&:= &\sum^{n-1}_{i=0}\frac{1}{t_{i+1}-t_i}
\int_{t_i}^{t_{i+1}}ds
\bigg \{\int_{t_i}^{s} \big  (\triangledown \sigma\cdot \sigma \big )(X_u(x))dB_u\bigg \}(B_{t_{i+1}}-B_{t_i})\\
&&- \frac{1}{2}\int_{0}^{t}\big ( \triangledown \sigma \cdot \sigma\big )(X_s(x))ds,\\
A_{3\pi}(x)&:=& \sum^{n-1}_{i=0}\frac{1}{t_{i+1}-t_i}
\int_{t_i}^{t_{i+1}}ds \bigg \{\int_{t_i}^{s} \big  (\triangledown \sigma\cdot {\bf \xi} \big )(X_u(x))
d|L^x|_u\bigg \}(B_{t_{i+1}}-B_{t_i})\\
A_{4\pi}(x)&:=& \sum^{n-1}_{i=0}\frac{1}{t_{i+1}-t_i}
\int_{t_i}^{t_{i+1}}ds \bigg \{\int_{t_i}^{s} \big  (\triangledown
\sigma\cdot b \big )(X_u(x))
du\bigg \}(B_{t_{i+1}}-B_{t_i})\\
A_{5\pi}(x)&:=& \frac{1}{2}\sum^{n-1}_{i=0}\frac{1}{t_{i+1}-t_i}
\int_{t_i}^{t_{i+1}}ds \bigg \{\int_{t_i}^{s} \big  (\triangledown
\sigma\cdot \triangledown \sigma \cdot \sigma \big )(X_u(x))
du\bigg \}(B_{t_{i+1}}-B_{t_i})\\
A_{6\pi}(x)&:=& \frac{1}{2}\sum^{n-1}_{i=0}\frac{1}{t_{i+1}-t_i}
\int_{t_i}^{t_{i+1}}ds \bigg \{\int_{t_i}^{s} {\bf tr}\big \{\big
(\triangledown^2 \sigma\cdot  \sigma \cdot \sigma^T \big )(X_u(x))
\big \}du\bigg \}(B_{t_{i+1}}-B_{t_i}).
\end{eqnarray*}
\begin{Prop}
Assume that the smooth bounded open $  \mathcal{O} $, the
coefficients $ \sigma $ and $b$ satisfy  the same conditions as in
Theorem 1.1.  $ ( X_t(x), L^x_t ) $ is a solution of Eq.(1.3). Then
for any $ p\geq 2$ and $R>0$ there exist  constant $c(p,R) $, which
is independent of $t$ and $ \pi $, and $\beta_0 \in (0,1)$
 such that
\begin{eqnarray}
\sup_{ x\in [-R,  R]^d\cap\bar{\mathcal{O}}}{\bf E} \big
\{|S_\pi(t,x)-I(t,x ) |^{2p}\big \}\leq c (p, R)\|\pi\|^{\beta_0 p}.
\end{eqnarray}
\end{Prop}
{\bf Proof.}  By Burkholder-Davis-Gundy  and H\"{o}lder
inequalities, we have
\begin{eqnarray*}
\{{\bf E} \{ |A_{1\pi}(x)|^{2p}\}\}^{\frac{1}{p}}&\leq & c(p) \bigg
[ {\bf E} \bigg (
\sum^{n-1}_{i=0}\int_{t_i}^{t_{i+1}}|(\sigma(X_s(x))-
\sigma(X_{t_i}(x))|^2ds\bigg )^p \bigg ]^{\frac{1}{p}  }\nonumber\\
&\leq & c(p)\sum^{n-1}_{i=0}\bigg \{ {\bf E}
\vert\int_{t_i}^{t_{i+1}}|(\sigma(X_s(x))-\sigma(X_{t_i}(x))|^2ds\vert^{p}
 \bigg
\}^{\frac{1}{p}}\nonumber\\
&\leq&c(p)\sum^{n-1}_{i=0}(t_{i+1}-t_i)^{1-\frac{1}{p}}\bigg [
\int_{t_i}^{t_{i+1}}{\bf E}|\sigma(X_s(x))
-\sigma(X_{t_i}(x))|^{2p}ds\bigg ]^{\frac{1}{p}}\nonumber \\
&\leq & c(p)\|\pi\|^{\frac{1}{2}},
\end{eqnarray*}
where we have used Proposition 2.2 and the condition (1.7). Thus
\begin{eqnarray}
{\bf E} \{ |A_{1\pi}(x)|^{2p}\}\leq c\|\pi\|^{\frac{p}{2}}.
\end{eqnarray}
Using Fubini Theorem,  $A_{2\pi}$ can be further written as
\begin{equation}
A_{2\pi}(x)=A_{2\pi}^{(1)}(x)+A_{2\pi}^{(2)}(x)+A_{2\pi}^{(3)}(x),
\end{equation}
where
\begin{eqnarray*}
A_{2\pi}^{(1)}(x)&:=& \sum^{n-1}_{i=0}\frac{1}{t_{i+1}-t_i}
\int_{t_i}^{t_{i+1}}(t_{i+1}-u) \big ( (\triangledown \sigma \cdot
\sigma )( X_u(x))- (\triangledown \sigma \cdot
\sigma )( X_{t_i}(x) )\big )du,\\
A_{2\pi}^{(2)}(x)&:=&
-\frac{1}{2}\sum^{n-1}_{i=0}\int_{t_i}^{t_{i+1}} \big (
(\triangledown \sigma \cdot \sigma )( X_u(x))- (\triangledown \sigma
\cdot \sigma )( X_{t_i}(x) )\big )du,\\
A_{2\pi}^{(3)}(x)&:= &\sum^{n-1}_{i=0}\bigg \{ \frac{1}{t_{i+1}-t_i}
\big (\int_{t_i}^{t_{i+1}}(t_{i+1}-u) (\triangledown \sigma \cdot
\sigma)(X_u(x))dB_u\big )(B_{t_{i+1}}-B_{t_i})\\
&&-\frac{1}{t_{i+1}-t_i} \big (\int_{t_i}^{t_{i+1}}(t_{i+1}-u)
(\triangledown \sigma \cdot \sigma )(X_u(x))du\big )\bigg
\}\nonumber\\
&:=&\sum^{n-1}_{i=0}( A_i +B_i).
\end{eqnarray*}
It follows from (1.7) and Proposition 2.2 that
\begin{eqnarray}
\{{\bf E} \{ |A_{2\pi}^{(1)}(x)|^{p}\}\}^{\frac{1}{p}}&\leq &
\sum^{n-1}_{i=0}\int_{t_i}^{t_{i+1}} \big \{{\bf E} \{\|\big (
(\triangledown \sigma \cdot \sigma )( X_u(x))- (\triangledown \sigma
\cdot \sigma )( X_{t_i}(x) )\big ) \|^p\}
\big \}^{\frac{1}{p}}du\nonumber\\
&\leq & c\sum^{n-1}_{i=0}\int_{t_i}^{t_{i+1}} \{{\bf E}
\{\| X_u(x)-X_{t_i}(x) \|^p\}\}^{\frac{1}{p}}du\nonumber\\
&\leq & c\|\pi\|^{\frac{1}{4}}.
\end{eqnarray}
Similar arguments lead to
\begin{equation}
\{{\bf E} \{ |A_{2\pi}^{(2)}(x)|^{p}\}\}^{\frac{1}{p}}\leq
c\|\pi\|^{\frac{1}{4}}.
\end{equation}
Since $A_{2\pi}^{(3)}$ is a  martingale  and $\triangledown \sigma
\cdot \sigma $ is bounded on $ \bar{\mathcal{O}} $, using
Burkholder-Davis-Gundy and H\"{o}lder inequalities, we obtain that
\begin{eqnarray}
&\{&{\bf E} \{ |A_{2\pi}^{(3)}(x)|^{2p}\}\}^{\frac{1}{p}}
\leq  c \bigg\{{\bf E} \big\{\sum^{n-1}_{i=0}| A_i +B_i|^2 \big\}^p
 \bigg\}^{\frac{1}{p}}\nonumber\\
&\leq & c \sum^{n-1}_{i=0} \big ( {\bf E}|A_i|^{2p} \big
)^{\frac{1}{p}}+c \sum^{n-1}_{i=0} \big ( {\bf E}|B_i|^{2p} \big
)^{\frac{1}{p}}\nonumber\\
&=& c(p)\sum^{n-1}_{i=0}\bigg \{ \big ({\bf E}\{ \big
|\frac{1}{t_{i+1}-t_i} \int_{t_i}^{t_{i+1}}
(t_{i+1}-u)(\triangledown \sigma \cdot \sigma )( X_u(x))dB_u\big
)\nonumber\\
&&\times(B_{t_{i+1}}
-B_{t_i})|^{2p}\big )^{\frac{1}{p}}\nonumber\\
&&+\big ({\bf E}\{ \big |\frac{1}{t_{i+1}-t_i} \big
(\int_{t_i}^{t_{i+1}}(t_{i+1}-u)(\triangledown \sigma \cdot \sigma )
( X_u(x))du\big |\}^{2p}\big )^{\frac{1}{p}}\bigg \}\nonumber\\
&\leq & c(p)\|\pi\|.
\end{eqnarray}
So we deduce from (4.7)-(4.9) that
\begin{eqnarray}
{\bf E} \{ |A_{2\pi}(x)|^{2p}\}\leq c\|\pi\|^{\frac{p}{4}}.
\end{eqnarray}
By Propositions 2.1 and 2.2, it follows from Kolmogorov's continuity
criterion( see Theorem 1.4.1 in \cite{s21} ) that there exist a
random variable $ K$ with $ {\bf E}|K(\omega)|^p < +\infty       $
and a positive constant $ \beta \in (0,1)$ such that
\begin{eqnarray}
\sup_{x\in [-R, R]^d \cap\bar{\mathcal{O}}}| L^x_{t}-L^x_s|\leq
K(\omega)|t-s|^\beta .
\end{eqnarray}
\begin{eqnarray}
{\bf E} \{ |A_{3\pi}(x)|^{2p}\} &\leq &{\bf E} \{ \big
(\sum^{n-1}_{i=0}\big |\frac{1}{t_{i+1}-t_i} \int_{t_i}^{t_{i+1}}ds
 \int_{t_i}^s(\triangledown \sigma \cdot {\bf \xi })( X_u(x) )d|L^x|_u\big|^2\big )^p\nonumber\\
&&\times \big (\sum^{n-1}_{i=0}|B_{t_{i+1}}-B_{t_i}|^2\big )^p\}\nonumber\\
&\leq & c{\bf E} \{ \big
(\sum^{n-1}_{i=0}(|L^x|_{t_{i+1}}-|L^x|_{t_i})^2\big )^p\times
 \big (\sum^{n-1}_{i=0}|B_{t_{i+1}}-B_{t_i}|^2\big )^p\}\nonumber\\
&\leq & c{\bf E} \{ \big (\sup_i|L_{t_{i+1}}^x-L_{t_i}^x|
|L^x|_1\big )^p\times
\big (\sum^{n-1}_{i=0}|B_{t_{i+1}}-B_{t_i}|^2\big )^p\}\nonumber\\
&\leq & c({\bf E}\{\big (\sup_i |L_{t_{i+1}}^x-L_{t_i}^x|\big
)^{3p}\})^{\frac{1}{3}}
 ({\bf E}\{\big |L^x\big |_{1}^{3p}\})^{\frac{1}{3}}\nonumber\\
&&\times \big ({\bf E} \{\big (\sum^{n-1}_{i=0}
|B_{t_{i+1}}-B_{t_i}|^2\big )^{3p}\}\big )^{\frac{1}{3}}\nonumber\\
&\leq& c \|\pi\|^{\beta p}
\end{eqnarray}
due to  the  inequality  (4.11) and  $\triangledown \sigma $ is
bounded on $ \bar{\mathcal{O}} $.
 For $p\geq 1$, by H\"{o}lder
inequality and  $ \triangledown^2 \sigma\cdot  \sigma \cdot \sigma^T
 $ is bounded on $ \bar{\mathcal{O}} $,
\begin{eqnarray}
&&{\bf E} \{ |A_{6\pi}(x)|^{2p}\}\nonumber\\
&\leq &{\bf E} \bigg\{ \big (\sum^{n-1}_{i=0}\big
|\frac{1}{t_{i+1}-t_i}
 \int_{t_i}^{t_{i+1}}ds \int_{t_i}^s{\bf Tr}\big \{\big
(\triangledown^2 \sigma\cdot  \sigma \cdot \sigma^T \big )(X_u(x))
\big \}du\big |^2\big )^p\nonumber\\
 &&\times \big (\sum^{n-1}_{i=0}|B_{t_{i+1}}-B_{t_i}|^2\big )^p\bigg\}\nonumber\\
&\leq & c\|\pi\|^p\sup_{\pi}{\bf E} \{\big
(\sum^{n-1}_{i=0}|B_{t_{i+1}}
-B_{t_i}|^2\big )^p\}\nonumber\\
&\leq& c\|\pi\|^p .
\end{eqnarray}
Similarly,
\begin{eqnarray}
{\bf E} \{ |A_{i\pi}(x)|^{2p}\}\leq c\|\pi\|^p, \ \mbox{ for $i=4,5$
}.
\end{eqnarray}
Putting the above estimates (4.5), (4.10) and (4.12)-(4.14) for
$A_{i\pi}$ ($i=1, \cdots, 6$) together, we deduce that
\begin{eqnarray*}
\sup_{ x\in [-R,  R]^d\cap\bar{\mathcal{O}}}{\bf E} \big
\{|S_\pi(t,x)-I(t,x ) |^{2p}\big \}\leq c (p, R)\|\pi\|^{\beta_0 p},
\end{eqnarray*}
where $ \beta_0=\min\{\frac{1}{4}, \beta  \}$. The proof of
Proposition 4.1 is complete. \ $\Box$ \vskip 0.3cm
 Next result is the moment estimates for
the two point motions.
\begin{Prop}
Assume that the smooth bounded open $  \mathcal{O} $, the
coefficients $ \sigma $ and $b$ satisfy  the same conditions as in
Theorem 1.1.  $ ( X_t(x), L^x_t ) $ is a solution of Eq.(1.3). Then
for any $ p\geq 2$ and $R>0$ there exists  constant $c(p,R) $, which
is independent of $t$ and $ \pi $,  such that
\begin{eqnarray}
{\bf E} \big \{\sup_{ t\in [0,  1]}|S_\pi(t,x)-S_\pi(t,y) |^{p}\big
\}\leq c(p, R ) |x-y|^{p},
\end{eqnarray}
for all $x, y\in [-R, R]^d\cap \bar{\mathcal{O}}$.
\end{Prop}
{\bf Proof.} Similarly as (4.3),
\begin{eqnarray}
S_\pi(t,x)-S_\pi(t,y) =\sum^6_{i=1}A_{i\pi}(x,y),
\end{eqnarray}
where
\begin{eqnarray*}
A_{1\pi}(x,y)&:=&\sum^{n-1}_{i=0}\big [ \sigma (X_{t_i}( x))- \sigma
(X_{t_i}( y))\big ]\big( B_{t_{i+1}}-B_{t_i}\big ),\\
A_{2\pi}(x,y )&:= &\sum^{n-1}_{i=0}\frac{1}{t_{i+1}-t_i}
\int_{t_i}^{t_{i+1}}ds \bigg \{\int_{t_i}^{s}\big [ \big
(\triangledown \sigma\cdot \sigma \big )(X_u(x))\\
&&-\big  (\triangledown \sigma\cdot \sigma \big )(X_u(y))\big ]dB_u
\bigg \}\big (B_{t_{i+1}}-B_{t_i}\big ),\\
A_{3\pi}(x,y )&:=& \sum^{n-1}_{i=0}\frac{1}{t_{i+1}-t_i}
\int_{t_i}^{t_{i+1}}ds \bigg \{\int_{t_i}^{s} \big  (\triangledown
\sigma\cdot {\bf \xi} \big )(X_u(x)) d|L^x|_u\\
&&-\int_{t_i}^{s} \big (\triangledown \sigma\cdot {\bf \xi} \big
)(X_u(y)) d|L^y|_u\bigg \}\big (B_{t_{i+1}}-B_{t_i}\big ),\\
A_{4\pi}(x,y)&:=& \sum^{n-1}_{i=0}\frac{1}{t_{i+1}-t_i}
\int_{t_i}^{t_{i+1}}ds \bigg \{\int_{t_i}^{s} \big [\big
(\triangledown \sigma\cdot b \big )(X_u(x))\\
&&-\big (\triangledown \sigma\cdot b \big )(X_u(y))\big ]du\bigg
\}\big (B_{t_{i+1}}-B_{t_i}\big ),
\end{eqnarray*}
\begin{eqnarray*}
 A_{5\pi}(x, y)&:=&
\frac{1}{2}\sum^{n-1}_{i=0}\frac{1}{t_{i+1}-t_i}
\int_{t_i}^{t_{i+1}}ds \bigg \{\int_{t_i}^{s} \big [\big
(\triangledown \sigma\cdot \triangledown \sigma \cdot \sigma \big
)(X_u(x))\\
&&- \big (\triangledown \sigma\cdot \triangledown \sigma \cdot
\sigma \big )(X_u(y))\big ]du\bigg \}\big (B_{t_{i+1}}-B_{t_i}\big ),\\
A_{6\pi}(x,y )&:=& \frac{1}{2}\sum^{n-1}_{i=0}\frac{1}{t_{i+1}-t_i}
\int_{t_i}^{t_{i+1}}ds \bigg \{\int_{t_i}^{s} {\bf Tr}\big \{\big
(\triangledown^2 \sigma\cdot  \sigma \cdot \sigma^T \big )(X_u(x))\\
&& -\big (\triangledown^2 \sigma\cdot  \sigma \cdot \sigma^T \big
)(X_u(y)) \big \}du\bigg \}\big(B_{t_{i+1}}-B_{t_i}\big ).
\end{eqnarray*}
Let $ | L^x- L^y |_t $ denote the total variation of $ L^x_\cdot-
L^y_\cdot$ on $ [0, t]$ for any $t\in [0,1]$, by (1.4), we have for
any $s\in [t_i, t_{i+1}]$
\begin{eqnarray*}
 &&\big |  \int_{t_i}^{s} \big  (\triangledown \sigma \cdot {\bf \xi}
\big )(X_u(x)) d|L^x|_u-\int_{t_i}^{s} \big (\triangledown
\sigma\cdot {\bf \xi} \big )(X_u(y)) d|L^y|_u    \big |^2\\
&=& \big |  \int_{t_i}^{s} \big  (\triangledown \sigma\big )(X_u(x))
dL^x_u-\int_{t_i}^{s} \big (\triangledown \sigma \big )(X_u(y))
dL^y_u \big |^2\\
&\leq & c \big |  \int_{t_i}^{s} \big [ (\triangledown \sigma
)(X_u(x)) -  (\triangledown \sigma  )(X_u(y))\big ]dL^x_u \big |^2+\big |  \int_{t_i}^{s} (\triangledown \sigma  )(X_u(y))d (  L^x_u-  L^y_u ) \big |^2\\
&\leq & c\sup_{s\in [0,1]}\{|X_s(x)-X_s(y) | \}\{ |L^x|_{t_{i+1}}-
|L^x|_{t_{i}} \} +c\{ |L^x-L^y|_{t_{i+1}}- |L^x-L^y|_{t_{i}} \}.
\end{eqnarray*}
So
\begin{eqnarray*}
(A_{3\pi}(x,y))^2&\leq &c\sum^{n-1}_{i=0}\frac{1}{t_{i+1}-t_i}
\int_{t_i}^{t_{i+1}}\big\{ |  \int_{t_i}^{s} \big  (\triangledown
\sigma \cdot {\bf \xi} \big )(X_u(x)) d|L^x|_u\nonumber\\
&& -\int_{t_i}^{s} \big (\triangledown \sigma\cdot {\bf \xi} \big
)(X_u(y)) d|L^y|_u ^2\big\}ds\times \sum^{n-1}_{i=0}|B_{t_{i+1}}-B_{t_i}|^2\nonumber\\
 &\leq & \bigg [ 2c^2 \sup_{s\in [0,1]} \{|X_s(x)-X_s(y)|^2\} |L^x|^2_1|
\nonumber\\
 &&+ 2c^2 \sup_{s\in [0,1]} \{|L^x_s-L^y_s|\} (|L^x|_1+ |L^y|_1)  \bigg
 ]\nonumber\\&&
 \times \sum^{n-1}_{i=0}|B_{t_{i+1}}-B_{t_i}|^2.
\end{eqnarray*}
It follows from Propositions 2.1-2.3 and H\"{o}lder's inequality
that
\begin{eqnarray}
&&{\bf E} \big \{\sup_{ t\in [0,  1]}|A_{3\pi}(x,y)
|^{2p}\big \}\nonumber\\
&\leq &c\big ( {\bf E} \big \{\sup_{0\leq s\leq 1} \vert
X_s(x)-X_s(y)\vert^{6p}\big \}\big )^{\frac{1}{3}}
\big ({\bf E} \big \{| L^x|_{1}^{6p}\big \})^{\frac{1}{3}}\nonumber\\
&&\times \big ({\bf E} \big \{\big (\sum^{n-1}_{i=0}
|B_{t_{i+1}}-B_{t_i}|^2\big )^{3p}\big
\}\big )^{\frac{1}{3}}\nonumber\\
&&+ c\big ( {\bf E} \big \{(\sup_{0\leq s\leq 1} |
L^x_s-L^y_s|^{3p}\big \}\big )^{\frac{1}{3}}
\big ({\bf E} \big \{| L^x|_{1}^{3p} +| L^y|_{1}^{3p} \big \})^{\frac{1}{3}}\nonumber\\
&&\times \big ({\bf E} \big \{\big (\sum^{n-1}_{i=0}
|B_{t_{i+1}}-B_{t_i}|^2\big )^{3p}\big
\}\big )^{\frac{1}{3}}\nonumber\\
&\leq & c |x-y|^{p}.
\end{eqnarray}
By Burkholder-Davis-Gundy inequalities, the condition (1.7) and
Proposition 2.1, it follows easily that
\begin{eqnarray}
{\bf E} \big \{\sup_{ t\in [0,  1]}|A_{1\pi}(x,y) |^{p}\big \}\leq
C(p, R ) |x-y|^{p}.
\end{eqnarray}
Using Burkholder-Davis-Gundy  and  H\"{o}lder's inequalities, the
condition (1.7) and Proposition 2.1,\\
\begin{eqnarray*}
\bigg ({\bf E} \big \{\sup_{ t\in [0,  1]}|A_{2\pi}(x,y) |^{p}\big
\}\bigg )^{\frac{1}{p}}\leq \sum^{n-1}_{i=0}\frac{1}{t_{i+1}-t_i}
\int_{t_i}^{t_{i+1}}ds \bigg \{{\bf E}\big\{\big |\int_{t_i}^{s}\big
[ \big
(\triangledown \sigma\cdot \sigma \big )(X_u(x))\nonumber\\
\end{eqnarray*}
\begin{eqnarray}
&&-\big  (\triangledown \sigma\cdot \sigma \big )(X_u(y))\big ]dB_u
\bigg \}\big (B_{t_{i+1}}-B_{t_i}\big )\big |^p\bigg \}^{\frac{1}{p}
}\nonumber\\
&\leq &\sum^{n-1}_{i=0}\big \{\frac{1}{t_{i+1}-t_i}
\int_{t_i}^{t_{i+1}}ds\bigg ( {\bf E} \big \{\vert
\int_{t_i}^{s}\big [(\triangledown \sigma\cdot \sigma \big )(X_u(x))
\nonumber\\
&&-(\triangledown \sigma\cdot \sigma \big )(X_u(y))\big
]dB_u\vert^{2p}\big \}\bigg )^{\frac{1}{2p}}\bigg ( {\bf E} \big
\{\vert (B_{t_{i+1}}
-B_{t_i})\vert^{2p}\big \}\bigg )^{\frac{1}{2p}}\big\}\nonumber\\
&\leq &\sum^{n-1}_{i=0}\big \{\frac{1}{t_{i+1}-t_i}
\int_{t_i}^{t_{i+1}}ds\bigg (\int_{t_i}^{s} \big \{{\bf E}[ \vert
X_u(x))-X_u(y)\vert^{2p}]
\big\}^{\frac{1}{p}}du\bigg )^{\frac{1}{2}}\nonumber\\
&&\times (t_{i+1}-t_i)^{\frac{1}{2}} \big\}\nonumber\\
 &\leq & c(p)|x-y|.
\end{eqnarray}
Similarly,
\begin{eqnarray}
&&\bigg ({\bf E} \big \{\sup_{ t\in [0,  1]}|A_{4\pi}(x,y)+
A_{5\pi}(x,y) |^{p}\big \}\bigg )^{\frac{1}{p}}\leq c(p)|x-y|,\\
&&\bigg ({\bf E} \big \{\sup_{ t\in [0,  1]}|A_{6\pi}(x,y) |^{p}\big
\}\bigg )^{\frac{1}{p}}\leq c(p)|x-y|.
\end{eqnarray}
 Thus, combining  above estimates (4.16)-(4.21) together , we  complete the
proof. \quad $\Box$ \vskip 0.3cm
 The following result can be proved similarly
as Proposition 4.2, but its proof is very easy, we omit it here.
\begin{Prop}
Assume that the smooth bounded open $  \mathcal{O} $, the
coefficients $ \sigma $ and $b$ satisfy  the same conditions as in
Theorem 1.1.  $ ( X_t(x), L^x_t ) $ is a solution of Eq.(1.3). Then
for any $ p\geq 2$ and $R>0$ there exists  constant $c(p,R) $, which
is independent of $t$ and $ \pi $,  such that
\begin{eqnarray}
{\bf E} \big \{\sup_{ t\in [0,  1]}|I(t,x)-I(t,y) |^{p}\big \}\leq
C(p, R ) |x-y|^{ p},
\end{eqnarray}
for all $x, y\in [-R, R]^d\cap \bar{\mathcal{O}}$.
\end{Prop}\vskip 0.3cm
\setcounter{equation}{0}
\section{{\bf{\small  Uniform convergence of the Riemann sums}}}
Let $R> 0$ and $p>1$ be given. Define $G_R:= [-R, R]^d\cap
\bar{\mathcal{O} }$. Then the following is a direct consequence of
Garsia-Rodemich and Rumsey's Lemma (cf.\cite{s24, s22,s13}).
\begin{Lemma}
Let $ f: \Omega \times \Re^m\rightarrow \Re^n $ be a measurable
stochastic field  taking values in $\Re^n$ which is continuous,
${\bf P}$- a.s., $p>1$. Then there exists a constant $c(p, R)$ such
that
\begin{eqnarray}
&&{\bf E}\big \{\sup_{x, y \in G_R}
\rho (f(x),f(y))^p\big\}\nonumber\\
&& \leq c(p, R)\int\int_{G_R\times G_R}{\bf E}\bigg\{ \frac{\rho
(f(x),f(y)|)^p }{d(x,y)^p}\bigg\}I_{\{x\neq y\}}dxdy,
\end{eqnarray}
where $(\Re^m, d )$ and $ ( \Re^n, \rho )$ are metric spaces.
\end{Lemma}
Now we prove the main result of this section.
\begin{theorem}
Assume that the smooth bounded open $  \mathcal{O} $, the
coefficients $ \sigma $ and $b$ satisfy  the same conditions as in
Theorem 1.1.  $ ( X_t(x), L^x_t ) $ is a solution of Eq.(1.3).
 Then
for any $ p\geq 2$ and $R>0$,
\begin{equation}
\lim_{\|\pi\|\rightarrow 0}{\bf E} \big \{\sup_{ x\in
G_R}|S_\pi(t,x)-I(t,x ) |^{2p}\big \}=0.
\end{equation}
\end{theorem}
{\bf Proof}. Since $ \sup_{x\in G_R}\{|f(x)|\}\leq \sup_{x, y\in
G_R}\{|f(x)-f(y)|\} + |f(x_0)|$ for any $x_0 \in G_R$ and function
$f$ on $ \Re^d$, we have
\begin{eqnarray}
&&{\bf E} \big \{\sup_{ x\in G_R}|S_\pi(t,x)-I(t,x ) |^{2p}\big \}\nonumber\\
 &&  \leq c(p){\bf E} \big \{\sup_{ x, y\in G_R}|S_\pi(t,x)-S_\pi(t,y)- I(t,x)
 +I(t,y)|^{2p}\big \} \nonumber\\
 &&+c(p) \sup_{ x\in G_R}{\bf E} \big \{|S_\pi(t,x)-I(t,x ) |^{2p}\big
 \}\nonumber\\
 &&\equiv  B_{1\pi }+ B_{2\pi} .
\end{eqnarray}
By Lemma 5.1,
\begin{eqnarray}
&& B_{1\pi }\nonumber\\
 && \leq c(p, R)\int\int_{G_R\times
G_R}{\bf E}\bigg\{ \frac{ |S_\pi(t,x)-S_\pi(t,y)- I(t,x)
+I(t,y)|^{2p} }{|x-y|^{2p}}\bigg\}I_{\{x\neq y\}}dxdy,\nonumber\\
\end{eqnarray}
By  Propositions 4.1,
\begin{eqnarray*}
 {\bf E}\bigg\{ \frac{ |S_\pi(t,x)-S_\pi(t,y)- I(t,x) +I(t,y)|^{2p}
}{|x-y|^{2p}}\bigg\}I_{\{x\neq y\}}\leq \frac{c\|\pi\|^{\beta_0
p}}{|x-y|^{2p}} I_{\{x\neq y\}}   \rightarrow 0
 \end{eqnarray*}
 as $\| \pi\|  \rightarrow 0$.\\   Therefore, by Propositions 4.2,
  dominated convergence theorem and (5.4), we have
\begin{eqnarray}
B_{1\pi }\longrightarrow 0\  , \mbox{ as} \ \pi\rightarrow 0.
\end{eqnarray}
By Proposition 4.1,
\begin{eqnarray}
B_{2\pi }\longrightarrow 0\  , \mbox{ as} \ \pi\rightarrow 0.
\end{eqnarray}
 Thus we complete the proof by (5.3),(5.5) and (5.6).
  \quad $\Box$ \vskip 0.3cm
\setcounter{equation}{0}
\section{ \bf{\small  Proof of Theorem 1.1} }
We will prove that $(X_t(Z), L_t^Z)$ solves the anticipating
reflected SDE (1.8).  Since $ ( X_t(x), L^x_t ) $ is a solution of
Eq.(1.3), by Proposition 3.1, we need only to prove (1.11). By
Theorem 5.1, the following  holds almost surely on $\{\omega;
Z(\omega)\in G_M\}$
$$ \int^t_0 \sigma (X_s(x))\circ dB_s\big |_{x=Z}\chi_{\{\omega; Z(\omega)\in G_M\}}$$
$$=\lim_{\|\pi\|\rightarrow 0}S_\pi(t,x)\big |_{x=Z}\chi_{\{\omega; Z(\omega)\in G_M\}}$$
$$= \int^t_0 \sigma (X_s(Z))\circ dB_s\chi_{\{\omega; Z(\omega)\in G_M\}}.$$
Letting $M\rightarrow \infty$, we obtain the substitution formula
(1.11), and therefore prove the Theorem. \quad \quad $\Box$
\begin{Remark}
By  checking carefully the proof of Theorem 1.1 and using Theorem
3.1 proved by Lions and Sznitman (see \cite{s3}),  the  conditions
on $\mathcal{O} $ in Theorem 1.1 can be weaken, that is, if
$\mathcal{O}$ satisfies the admissibility condition (see \cite{s3},
page 521) and the following condition: there exists a function $\phi
$ in $\mathcal{C}^2_b(\Re^d)$ such that $  \exists \alpha >0 $,
$\forall x \in \partial \mathcal{O} $, $  \forall y \in
\bar{\mathcal{O} }$, $\forall \xi \in {\bf n}(x)$\  $ \Rightarrow
\frac{1}{\alpha} (\triangledown \phi (x), \xi )|y-x|^2-(y-x, \xi
)\leq 0 $,  Theorem 1.1 also holds.
\end{Remark}
 \vskip 1cm
 {\bf Acknowledgements.} This work is supported by NSFC and SRF
for ROCS, SEM. The author would like to thank both for their
generous financial support.

\end{document}